\documentclass[a4paper,10pt]{article}
\usepackage[utf8]{inputenc}
\usepackage{amssymb}
\usepackage{amsthm}
\usepackage{tikz}
\usepackage{tikz-cd}
\usepackage{bm} 
\usepackage{color}
\usepackage{mathtools}
\usepackage[titletoc,page]{appendix}
\usepackage[export]{adjustbox}
\usepackage{titletoc}
\usepackage{enumitem}

\usepackage[mathscr]{euscript} 
\DeclareFontFamily{OT1}{pzc}{}
\DeclareFontShape{OT1}{pzc}{m}{it}{<-> s * [1.10] pzcmi7t}{}
\DeclareMathAlphabet{\mathpzc}{OT1}{pzc}{m}{it} 

\dottedcontents{section}[2.5em]{\bfseries}{2.5em}{1pc} 

\usetikzlibrary{calc}
\usetikzlibrary{matrix,arrows,decorations}
\newtheorem{theorem}{Theorem}[section]
\newtheorem{lemma}[theorem]{Lemma}
\newtheorem{remark}[theorem]{Remark}

\newtheorem{example/theorem}[theorem]{Example/Theorem}
\newtheorem{construction/theorem}[theorem]{Construction/Theorem}

\newtheorem{convention/notation}[theorem]{Convention/Notation}
\newtheorem{corollary}[theorem]{Corollary}
 
\newtheorem{construction/definition}[theorem]{Construction/Definition}  
\newtheorem{notation}[theorem]{Notation}

\newtheorem{conjecture}[theorem]{Conjecture}

\usepackage{hyperref}
\title{An approximation of the $e$-invariant in the stable homotopy category}  
\author{Yi-Sheng Wang\footnote{yswangbl@gmail.com}}  

\begin{document}
\maketitle
\begin{abstract}  
In their construction of the topological index for flat vector bundles, Atiyah, Patodi and Singer associate to each flat vector bundle a particular $\mathbb{C/Z}$-$K$-theory class. This assignment determines a map, up to weak homotopy, from $K_{a}\mathbb{C}$, the $0$-connective algebraic $K$-theory space of the complex numbers, to $F_{t,\mathbb{C/Z}}$, the homotopy fiber of the Chern character. In this paper we give evidence for the conjecture that this map can be represented by an infinite loop map. The result of the paper implies a refined Bismut-Lott index theorem for a compact smooth bundle $E\rightarrow B$ with the fundamental group $\pi_{1}(E,\ast)$ finite, for every point $\ast\in E$. This paper is a continuation of the author's paper ``Topological $K$-theory with coefficients and the $e$-invariant".
\end{abstract}

\tableofcontents

\section{Introduction}
Throughout the paper, we use the Quillen model category of pointed $k$-spaces $\mathpzc{Top}_{\ast}$ as our convenient category of topological spaces. Given $X, Y\in\mathpzc{Top}_{\ast}$, $[X,Y]$ denotes the homotopy classes of maps from $X$ to $Y$ in $\mathpzc{Top}_{\ast}$.

In their construction of the topological index for flat vector bundles, Atiyah et al. associate to each flat vector bundle a $\mathbb{C/Z}$-$K$-theory class \cite[p.89]{APS3}, and this assignment gives a homomorphism $\bar{e}_{\operatorname{APS}}:\tilde{K}(M,\mathbb{C})\rightarrow [M,F_{t,\mathbb{C/Z}}]$, where $\tilde{K}(M,\mathbb{C})$ is the abelian group of zero dimensional virtual flat vector bundles over $M$, a compact smooth manifold \cite[Remark $4.1.3$]{Wang1}. Via the universal property of the plus construction, we further obtain a map $e:K_{a}\mathbb{C}\rightarrow F_{t,\mathbb{C/Z}}$ \cite{JW} and \cite{Wang1}. The construction can be summarized in the commutative diagram below:  
\begin{center}
\begin{equation}\label{Intro:Diag1}
\begin{tikzpicture}[baseline=(current  bounding  box.center)] 
\node(Lu) at (0,1.5) {$\tilde{K}(M,\mathbb{C})$};
\node(Ru) at (6,1.5) {$[M,F_{t,\mathbb{C/Z}}]$}; 
\node(Ml) at (3,0) {$[M,K_{a}\mathbb{C}]$};
\draw [->] (Lu) to [out=-90,in=180]  (Ml);
\draw [->] (Ml) to [out= 0,in = -90] node [yshift=.5em] {\scriptsize $e_{\ast}$}(Ru);   
\draw [->] (Lu) to node [above]{\scriptsize $\bar{e}_{\operatorname{APS}}$}(Ru);     
\end{tikzpicture}
\end{equation}
\end{center} 

It is also proved in \cite[Theorem $3.1$]{JW} (see \cite[Theorem $4.3.5$]{Wang1}) that the Borel regulator $\operatorname{Bo}$, a linear combination of Borel classes, factors through the map $e$ as follows: 
\begin{center}
\begin{tikzpicture} 
\node(L)  at (0,2) {$K_{a}\mathbb{C}$};
\node(M)  at (3,2) {$F_{t,\mathbb{C/Z}}$};
\node(MR) at (5,2) {$\Omega K_{t,\mathbb{R}}$};
\node(R)  at (7,2) {$H^{odd}\mathbb{R}$};
 
\draw [->]  (L) to node [above] {\scriptsize $e$}(M);  
\draw [->]  (M) to node [above] {\scriptsize $\operatorname{Im}$}(MR);
\draw [->]  (MR)to node [above] {\scriptsize $\operatorname{ch}_{\otimes\mathbb{R}}$} (R); 
\draw [->]  (L) to [out=30,in=150] node [above]{\scriptsize $\operatorname{Bo}$}(R); 
\end{tikzpicture}
\end{center}
where $H^{odd}\mathbb{R}:=\prod\limits_{\mathclap{i \text{ odd }}}K(\mathbb{R},i)$ and $K(A,i)$ denotes the Eilenberg-Maclane space of an abelian group $A$ in degree $i$. The maps $\operatorname{Im}:F_{t,\mathbb{C/Z}}\rightarrow \Omega K_{t,\mathbb{R}}$ and $\operatorname{ch}_{\otimes \mathbb{R}}:\Omega K_{t,\mathbb{R}}\xrightarrow{\sim} H^{odd}\mathbb{R}$
are induced by the homomorphism
\begin{align*}
\mathbb{C/Z}&\rightarrow \mathbb{R}\\
a+ib&\mapsto b
\end{align*}
and the Chern character, respectively, where $K_{t,\mathbb{R}}$ is the infinite loop space representing $0$-connective complex topological $K$-theory with coefficients in $\mathbb{R}$. 
Now the Bismut-Lott index theorem \cite{BL} says, for every compact smooth fiber bundle $E\rightarrow B$, the following diagram of abelian groups commutes: 
\begin{center}
\begin{tikzpicture} 
\node(Lu) at (0,2) {$\tilde{K}(E,\mathbb{C})$};
\node(Ll) at (0,0) {$\tilde{K}(B,\mathbb{C})$}; 
\node(Ru) at (4,2) {$[E,H^{odd}\mathbb{R}]$};
\node(Rl) at (4,0) {$[B,H^{odd}\mathbb{R}]$};

\draw[->] (Lu) to node [right]{\scriptsize $\pi^{!}$}(Ll);
\draw[->] (Ru) to node [right]{\scriptsize $\operatorname{tr}_{\operatorname{BG}}^{\ast}$}(Rl);
\draw[->] (Lu) to node [above]{\scriptsize $\bar{\operatorname{Bo}}$}(Ru);
\draw[->] (Ll) to node [above]{\scriptsize $\bar{\operatorname{Bo}}$}(Rl);  
\end{tikzpicture}
\end{center}
where $\pi^{!}$ is given by taking the fiberwise homology of $E\rightarrow B$ with coefficients in a flat vector bundle, $\operatorname{tr}_{\operatorname{BG}}$ is the Beck-Gottlieb transfer and $\bar{\operatorname{Bo}}$ is the composition
\[\tilde{K}(-,\mathbb{C})\rightarrow [-,K_{a}\mathbb{C}]\xrightarrow{\operatorname{Bo}_{\ast}} [-,H^{odd}\mathbb{R}].\]

We conjecture that there should be a finer index theorem in terms of $\bar{e}_{\operatorname{APS}}$:
\begin{conjecture}\label{Intro:Conj1}
The following diagram is commutative
\begin{center}
\begin{tikzpicture} 
\node(Lu) at (0,2) {$\tilde{K}(E,\mathbb{C})$};
\node(Ll) at (0,0) {$\tilde{K}(B,\mathbb{C})$}; 
\node(Ru) at (4,2) {$[E,F_{t,\mathbb{C/Z}}]$};
\node(Rl) at (4,0) {$[B,F_{t,\mathbb{C/Z}}]$};

\draw[->] (Lu) to node [right]{\scriptsize $\pi^{!}$}(Ll);
\draw[->] (Ru) to node [right]{\scriptsize $\operatorname{tr}_{\operatorname{BG}}^{\ast}$}(Rl);
\draw[->] (Lu) to node [above]{\scriptsize $\bar{e}_{\operatorname{APS}}$}(Ru);
\draw[->] (Ll) to node [above]{\scriptsize $\bar{e}_{\operatorname{APS}}$}(Rl);  
\end{tikzpicture}
\end{center}
\end{conjecture}
\noindent
This implies the $\operatorname{BL}$ index theorem after composing $\operatorname{ch}_{\otimes\mathbb{R},\ast}\circ J_{\ast}$.

In view of the diagram \eqref{Intro:Diag1} and the Dwyer-Weiss-Williams index theorem \cite{DWW}, which entails commutativity of the diagram below:
\begin{center}
\begin{equation}\label{Intro:DWWindex}
\begin{tikzpicture}[baseline=(current bounding box.center)]
\node(Lu) at (0,2) {$\tilde{K}(E,\mathbb{C})$};
\node(Ll) at (0,0) {$\tilde{K}(B,\mathbb{C})$};
 
\node(Ru) at (4,2) {$[E,K_{a}(\mathbb{C})]$};
\node(Rl) at (4,0) {$[B,K_{a}(\mathbb{C})]$};

\path[->, font=\scriptsize,>=angle 90]  
 
(Lu) edge (Ru)
(Ll) edge (Rl)    
(Lu) edge node [right]{$\pi^{!}$}(Ll)
(Ru) edge node [right]{$\operatorname{tr}_{\operatorname{BG}}^{\ast}$}(Rl); 
\end{tikzpicture}
\end{equation}
\end{center}
Conjecture \ref{Intro:Conj1} ensues from the following conjecture:
\begin{conjecture}\label{Intro:Conj2}
The map $e$ is weakly homotopic to an infinite loop map.
\end{conjecture}
\noindent
The results of this paper give evidence in support of Conjecture \ref{Intro:Conj2}.
 
\subsection*{Main results} 
\begin{theorem}\phantomsection\label{Intro:enatural}
Let $K_{a}\mathbb{C}\rightarrow K_{t}$ be the (canonical) comparison map from the $0$-connective algebraic $K$-theory space of the complex numbers to the $0$-connective complex topological $K$-theory space. Then there exists an infinite loop map $e^{\natural}_{h}$, unique up to phantom maps, such that the composition 
\[K_{a}\mathbb{C}\xrightarrow{e_{h}^{\natural}} F_{t,\mathbb{C/Z}}\rightarrow K_{t}\]
and the (canonical) comparison map $K_{a}\mathbb{C}\rightarrow K_{t}$ are homotopic as infinite loop maps and the map $e^{\natural}_{h}$ satisfies
\begin{align*}
e^{\natural}_{h,\ast}=&e_{\ast}:\pi_{\ast}(K_{a}\mathbb{C})\rightarrow \pi_{\ast}(F_{t,\mathbb{C/Z}}),\nonumber\\
e^{\natural}_{h,\ast}\vert_{\operatorname{Tor}}=&e_{\ast}\vert_{\operatorname{Tor}}:\operatorname{Tor}[L,K_{a}\mathbb{C}]\rightarrow [L,F_{t,\mathbb{C/Z}}],\\
\operatorname{ch}_{\otimes\mathbb{R}}\circ \operatorname{Im}\circ e^{\natural}_{h}=&\operatorname{Bo}\in \operatorname{Ho}(\mathcal{P}),
\end{align*} 
for every finite $CW$-complex $L$. 
\end{theorem}
\noindent 
Theorem \ref{Intro:enatural} shows that, in a certain sense, the map $e_{h}^{\natural}$ is the unique approximation of the map $e$ in the stable homotopy category, and it also implies the following index theorem, which is slightly weaker than Conjecture \ref{Intro:Conj1}: 
\begin{theorem}\label{Intro:Indexthm}
Let $E\rightarrow B$ be a smooth compact fiber bundle. If the fundamental group $\pi_{1}(E,\ast)$ is finite, for any base point $\ast\in E$, then the diagram    
\begin{center}
\begin{equation} 
\begin{tikzpicture}[baseline=(current  bounding  box.center)]
\node(Lu) at (0,2) {$\tilde{K}(E,\mathbb{C})$};
\node(Ll) at (0,0) {$\tilde{K}(B,\mathbb{C})$}; 
\node(Ru) at (4,2) {$[E,F_{t,\mathbb{C/Z}}]$};
\node(Rl) at (4,0) {$[B,F_{t,\mathbb{C/Z}}]$};

\draw[->] (Lu) to node [right]{\scriptsize $\pi^{!}$}(Ll);
\draw[->] (Ru) to node [right]{\scriptsize $\operatorname{tr}_{\operatorname{BG}}^{\ast}$}(Rl);
\draw[->] (Lu) to node [above]{\scriptsize $\bar{e}_{\operatorname{APS}}$}(Ru);
\draw[->] (Ll) to node [above]{\scriptsize $\bar{e}_{\operatorname{APS}}$}(Rl);  
\end{tikzpicture}
\end{equation}
\end{center} 
commutes
\end{theorem}

The next theorem shows that the map $e$ can be viewed as a generalization of the Adams $e$-invariant.\footnote{This relation between the map $e$ and the Adams $e$-invariant has been claimed without proof in \cite[p.930]{JW}. The name $e$ is also due to Jones and Westbury.}  
\begin{theorem}\label{Intro:eandeAdam}
Let $e_{\operatorname{Adams},\ast}:\pi_{\ast}^{s}=\pi_{\ast}(B\Sigma_{\infty}^{+})\rightarrow \mathbb{Q/Z}\subset\mathbb{C/Z}$ be the Adams $e$-invariant. Then we have
\[e_{\ast}\circ \iota_{\ast}=e_{\operatorname{Adams},\ast},\] 
where $B\Sigma_{\infty}^{+}$ is the plus construction of the classifying space of the infinite symmetric group $\Sigma_{\infty}$ and $\iota$ is the map induced from the canonical embedding $\Sigma_{\infty}\rightarrow \operatorname{GL}(\mathbb{C})$, where $\operatorname{GL}(\mathbb{C})$ is the infinite general linear group.  
\end{theorem} 
\subsection*{Outline of the paper}
To approximate the map $e:K_{a}\mathbb{C}\rightarrow F_{t,\mathbb{C/Z}}$ by infinite loop maps, we consider liftings---dashed arrow---of the following diagram: 
\begin{center}
\begin{tikzpicture}
\node (Lm) at (0,0){$K_{a}\mathbb{C}$};
\node (Ru) at (3,1){$F_{t,\mathbb{C/Z}}$};
\node (Rm) at (3,0){$K_{t}$};
\node (Rl) at (3,-1){$H^{ev}\mathbb{C}$};

\draw [->] (Lm) to (Rm) ;
\draw [dashed,->] (Lm) to node [above]{\scriptsize $e_{h}^{\mathbb{C}}$}(Ru);
\draw [->] (Ru)--(Rm);
\draw [->] (Rm) to node [right]{\scriptsize $\operatorname{ch}$}(Rl);

\end{tikzpicture}
\end{center}
where $K_{a}\mathbb{C}\rightarrow K_{t}$ is the (canonical) comparison map, $\operatorname{ch}$ is the Chern character, $H^{ev}\mathbb{C}:=\prod\limits_{\mathclap{i \text{ even }}}K(\mathbb{C},i)$, and the sequence 
\[F_{t,\mathbb{C/Z}}\rightarrow K_{t}\xrightarrow{\operatorname{ch}} H^{ev}\mathbb{C}\]
is a homotopy fiber sequence. In Section $1$, we shall see the existence of liftings $e^{\mathbb{C}}_{h}$ in the stable homotopy category and how their induced homomorphisms all restrict to $e_{\ast}$ on the torsion subgroup of $[L,K_{a}\mathbb{C}]$, for every finite $CW$-complex $L$. We show there are infinitely many different liftings in the stable homotopy category in Section $2$. In Section $3$, we investigate the relation between liftings $e_{h}^{\mathbb{C}}$ and their induced maps $t_{h}^{\mathbb{C}}:K^{rel}\mathbb{C}\rightarrow  H^{odd}\mathbb{C}$, where $K^{rel}\mathbb{C}$ is the relative algebraic $K$-theory space of the complex numbers, the homotopy fiber of the comparison map $K_{a}\mathbb{C}\rightarrow K_{t}$, and $H^{odd}\mathbb{C}:=\prod\limits_{\mathclap{i \text{ odd }}}K(\mathbb{C},i)$. Utilizing the relation between $e^{\mathbb{C}}_{h}$ and $t^{\mathbb{C}}_{h}$, we construct the map $e^{\natural}_{h}$ and prove the main theorem (\ref{Intro:enatural}). Finally, we apply the results obtained in the previous sections to prove the index theorem \eqref{Intro:Indexthm} and the comparison theorem \eqref{Intro:eandeAdam} in Sections $4$ and $5$. 
 
\subsection*{Relation and comparison to Bunke's regulators}
In \cite[Section $2.5$]{Bu1} and \cite[Definition $13.13$]{Bu2}, two maps from $K_{a}\mathbb{C}$ to $F_{t,\mathbb{C/Z}}$, which the comparison map $K_{a}\mathbb{C}\rightarrow K_{t}$ factors over, are constructed directly in the stable homotopy category via the technique of the $\infty$-categorical approach to $K$-theory developed in \cite{BNV} and \cite{BT}. Notice that their method actually give maps between $(-1)$-connective spectra, but one can always lift maps between $(-1)$-connective spectra to maps between (corresponding) $0$-connective spectra. The constructions of these two maps are very concrete and can be generalized to more general settings. They both involve geometries of vector bundles and are closely linked to differential $K$-theory. In fact, it is the use of the geometries, along with the $\infty$-categorical technique, gives the preferred homotopies needed to define the maps from $K_{a}\mathbb{C}$ to $F_{t,\mathbb{C/Z}}$. Though the constructions of these two maps and the map $e$ are similar in many respects, the precise relation between them is not entirely clear, and one might need a space level comparison in order to unravel it (see \cite[Remark $12.15$]{Bu2}). 

On the other hand, in this paper, we use a pure homotopy-theoretic method to find the approximation of the map $e$ in the stable homotopy category. Since there is no geometries of vector bundles involved in this approach, it is not easy to find a preferred homotopy. Therefore, instead of searching for a preferred homotopy, we study all the maps from $K_{a}\mathbb{C}$ to $F_{t,\mathbb{C/Z}}$ over which the comparison map $K_{a}\mathbb{C}\rightarrow K_{t}$ factors---there are infinitely many such maps (Theorem \ref{Manyliftings}). The idea is to find the homotopy properties that distinguish these maps from each other (up to weak homotopy). It proceeds roughly as follows: We first show all such maps restrict to the same map on the torsion part of $K_{a}\mathbb{C}$ (Theorem \ref{sameisoLiftings} and Corollary \ref{Torsionpartarethesame}) and hence what determines each of them is its restriction to the rational part of $K_{a}\mathbb{C}$. This observation, in view of the structure theorem of (relative) algebraic $K$-theory (Corollaries \ref{StrKa} and \ref{StrKr}), leads us to consider relative $K$-theory of the complex numbers and the approximation $t_{h}^{\natural}$ of the Chern-Simons class associated to the Chern character in the stable homotopy category. Note that the Chern-Simons classes give a map from $K^{rel}\mathbb{C}$ to $H^{odd}\mathbb{C}$ (see \cite[Section $4.3$]{Wang1}), which, when restricted to the rational part of $K_{a}\mathbb{C}$, can be considered as a (weak) lifting of the map $e$ and the map $t_{h}^{\natural}$ is an approximation of the Chern-Simons class in the sense that it induces the same homomorphism that Chern-Simons class induces between the homotopy groups of $K^{rel}\mathbb{C}$ and $H^{odd}\mathbb{C}$.
Then, using the connection between $t_{h}^{\mathbb{C}}:K^{rel}\mathbb{C}\rightarrow H^{odd}\mathbb{C}$ and $e_{h}^{\mathbb{C}}:K_{a}\mathbb{C}\rightarrow F_{t,\mathbb{C/Z}}$ (the diagram \eqref{tecommutativediag}, Corollary \ref{edeterminest} and Lemma \ref{thdetermineseh}), we find the approximation $e_{h}^{\natural}$ of the map $e$ and prove its uniqueness (Theorem \ref{enatural}). 

In particular, if one can show either of the maps constructed in \cite[Section $2.5$]{Bu1} and \cite[Definition $13.13$]{Bu2} is weakly homotopic to the map $e$, then the maps $e_{h}^{\natural}$ and $e$ are weakly homotopic as well. The latter can also be deduced if the map $t_{h}^{\natural}$ indeed realizes the Chern-Simons class, for every compact smooth manifold (Lemma \ref{thdetermineseh}).

\subsection*{Notation and conventions}
In this paper, we use the homotopy category of prespectra $\operatorname{Ho}(\mathcal{P})$ (see \cite[Appendix]{Wang1}) as our model for the stable homotopy category, and as in \cite{Wang1}, bold letters are reserved for prespectra and maps between them. Given two infinite loop spaces $E_{0}=\Omega^{\infty}\mathbf{E}$ and $F_{0}=\Omega^{\infty} \mathbf{F}$, where $\mathbf{E}, \mathbf{F}\in \mathcal{P}$, we let $[E_{0},F_{0}]_{\operatorname{Ho}(\mathcal{P})}:=[\mathbf{E},\mathbf{F}]_{\operatorname{Ho}(\mathcal{P})}$, the abelian group of maps between cofibrant-fibrant replacements of $\mathbf{E}$ and $\mathbf{F}$. When we say a diagram of infinite loop spaces is commutative, it commutes in $\operatorname{Ho}(\mathcal{P})$, unless otherwise specified.

Since, in most cases, our methods work in a more general setting that includes algebraic $K$-theory of the real numbers, real topological $K$-theory and topological $K$-theory with coefficients in $\mathbb{Q/Z}$ or $\mathbb{R/Z}$, we introduce the following notation for the sake of convenience:
\begin{notation}\phantomsection\label{notationsforliftings} 
\begin{itemize}[itemsep=1.5pt]
\item[] $\mathbb{F}^{\prime}=\mathbb{R}$ or $\mathbb{C}$, and $\mathbb{F}=\mathbb{Q}$ or $\mathbb{F}^{\prime}$.
\item[$K_{a}\mathbb{F}^{\prime}$:] The infinite loop space of the $\Omega$-$CW$-prespectrum $\mathbf{K}_{a}\mathbb{F}^{\prime}$ that represents $0$-connective algebraic $K$-theory of the complex numbers (resp. the real numbers). We often drop the field $\mathbb{F}^{\prime}$ from the notation when both cases apply. 

\item[$K_{t}$:] The infinite loop space of the $\Omega$-$CW$-prespectrum $\mathbf{K}_{t}$ that represents $0$-connective complex (resp. real) topological $K$-theory.  

\item[] When a statement is true for both $\mathbb{R}$ and $\mathbb{C}$, we do not specify the field. For instance, the map $K_{a}\rightarrow K_{t}$ could mean the comparison map from the $0$-connective algebraic $K$-theory space of the complex numbers to the $0$-connective complex topological $K$-theory space or the comparison map from the $0$-connective algebraic $K$-theory space of the real numbers to the $0$-connective real topological $K$-theory space. If a statement applies to just one case, we specify only the field used in algebraic $K$-theory. For example, the map $K_{a}\mathbb{C}\rightarrow K_{t}$ stands for the comparison map from the $0$-connective algebraic $K$-theory space of the complex numbers to the $0$-connective complex topological $K$-theory space.

\item[$X_{\mathbb{F}}$:] The infinite loop space of the prespectrum $\mathbf{X}\wedge \mathbf{M}\mathbb{F}$, or equivalently the zero component of its fibrant replacement, where $\mathbf{X}$ is a $CW$-prespectrum. In the case where $\mathbf{X}=\mathbf{K}_{t}$, we have the homotopy equivalence $K_{t,\mathbb{F}}\simeq H^{ev}\mathbb{F}$ given by the Chern character.
 
\item[$F_{t,\mathbb{F}/\mathbb{Z}}$:] The homotopy fiber of $K_{t}\rightarrow K_{t,\mathbb{F}}$, or equivalently the infinite loop space of the prespectrum $\Omega(\mathbf{K}_{t}\wedge \mathbf{M}\mathbb{F}/\mathbb{Z})$ (see \cite[Lemma $A.1.3$-$4$]{Wang1}).

\item[$F_{a,\mathbb{Q}/\mathbb{Z}}$:] The homotopy fiber of $K_{a}\rightarrow K_{a,\mathbb{Q}}$, or equivalently the infinite loop space of the prespectrum $\Omega(\mathbf{K}_{a}\wedge \mathbf{M}\mathbb{Q}/\mathbb{Z})$ (see \cite[Lemma $A.1.3$-$4$]{Wang1}).

\item[$K^{rel}\mathbb{F}^{\prime}$:] The homotopy fiber of $K_{a}\mathbb{F}^{\prime}\rightarrow K_{t}$. It is the infinite loop space representing relative $K$-theory of $\mathbb{F}^{\prime}$. The field $\mathbb{F}^{\prime}$ is dropped from the notation when a statement holds for both $\mathbb{R}$ and $\mathbb{C}$.  
\end{itemize}
\end{notation}

%%%Working in the stable homotopy category 

%%%%%%%%%%%Be careful about one thing: How can you be sure the Suslin's map is the one induced by the cofiber construction.

%%%%%%%%%%%%%%%%%%%%%%%%%%%%%  

%Beside there is a similar lemma as in the case of the Moore spaces:

%\begin{lemma}
%For any homomorphism $G\rightarrow \pi_{n}(\mathbf{E})$, one can find a map of prespectra 
%\[\mathbf{M}G[n] \rightarrow \mathbf{E}\] 
%such that it realizes the homomorphism after applying $\pi_{n}(-)$.
%\end{lemma}
%\begin{proof}
%From the standard construction or see \cite[4.33]{Ru}
%\end{proof}
%With this, the uniqueness of the homotopy type of a Moore prespectrum can be verified.

\section{Liftings of the comparison map $K_{a}\rightarrow K_{t}$}  
The comparison map $\mathbf{K}_{a}\rightarrow \mathbf{K}_{t}$ along with the generalized Chern character $\mathbf{E}\wedge \mathbf{M}\mathbb{Z}\rightarrow \mathbf{E}\wedge \mathbf{M}\mathbb{F}$ induced by the inclusion $\mathbb{Z}\hookrightarrow \mathbb{F}$, where $\mathbf{E}\in\mathcal{P}$, gives us the following commutative diagram of prespectra. 
%%%%Some issue about how to define K_{t}\rightarrow K_{t}MQ
\begin{center}
\begin{tikzpicture}
 \node (AK-1) at (0,8) {$\Omega(\mathbf{K}_{a}\wedge \mathbf{M}\mathbb{Z})$};
 \node (AQ-1) at (0,6) {$\Omega(\mathbf{K}_{a}\wedge\mathbf{M}\mathbb{Q})$};
 \node (FAQ)  at (0,4) {$\mathbf{Fib}(f)$};
 \node (AK)   at (0,2) {$\mathbf{K}_{a}\wedge  \mathbf{M}\mathbb{Z}$};
 \node (AQ)   at (0,0) {$\mathbf{K}_{a}\wedge \mathbf{M}\mathbb{Q}$};

 \node (TK1-1) at (4,8){$\Omega(\mathbf{K}_{t}\wedge  \mathbf{M}\mathbb{Z})$};
 \node (TQ-1) at (4,6) {$\Omega(\mathbf{K}_{t} \wedge  \mathbf{M}\mathbb{Q})$};
 \node (FTQ)  at (4,4) {$\mathbf{Fib}(g)$};
 \node (TK1)  at (4,2) {$\mathbf{K}_{t}\wedge \mathbf{M}\mathbb{Z}$};
 \node (TQ)   at (4,0) {$\mathbf{K}_{t} \wedge  \mathbf{M}\mathbb{Q}$};

 \node (TK2-1)at (8,8) {$\Omega(\mathbf{K}_{t}\wedge  \mathbf{M}\mathbb{Z})$};
 \node (TF-1) at (8,6) {$\Omega(\mathbf{K}_{t} \wedge  \mathbf{M}\mathbb{F}^{\prime})$};
 \node (FT)   at (8,4) {$\mathbf{Fib}(h)$};
 \node (TK2)  at (8,2) {$\mathbf{K}_{t}\wedge  \mathbf{M}\mathbb{Z}$};
 \node (TF)   at (8,0) {$\mathbf{K}_{t} \wedge  \mathbf{M}\mathbb{F}^{\prime}$};

 \path [->, font=\scriptsize, >=angle 90]
 (AK-1) edge node [right]{$\Omega \mathbf{f}$}(AQ-1) 
 (AQ-1) edge (FAQ) 
 (FAQ)  edge (AK)
 (AK)   edge node [right]{$\mathbf{f}$} (AQ)
 
 (TK1-1) edge node [right]{$\Omega \mathbf{g}$}(TQ-1)
 (TQ-1) edge (FTQ)
 (FTQ)  edge (TK1)
 (TK1)  edge node [right]{$\mathbf{g}$}(TQ)
 
 (TK2-1)  edge node [right]{$\Omega \mathbf{h}$}(TF-1)
 (TF-1) edge (FT)
 (FT)   edge (TK2)
 (TK2)  edge node [right]{$\mathbf{h}$}(TF)

 (AK-1) edge (TK1-1)
 (TK1-1) edge (TK2-1)
 (AQ-1) edge (TQ-1)
 (TQ-1) edge (TF-1) 
 (FAQ)  edge node [above]{$\mathbf{\operatorname{Su}}$} node[below]{$\sim$} (FTQ)
 (FTQ)  edge (FT) 
 (AK)   edge (TK1)
 (TK1)  edge node [above]{$id$}(TK2)
 (AQ)   edge (TQ)
 (TQ)   edge (TF) ; 
 
\end{tikzpicture}
\end{center}
where $\mathbf{\operatorname{Su}}$ stands for the $\pi_{\ast}$-isomorphism in \cite[Corollary $2.1.6$]{Wang1}. Using the model structure of $\mathcal{P}$ \cite[Theorem $A.1.2$]{Wang1}, one can replace each $CW$-prespectrum in the diagram above with an equivalent fibrant-cofibrant prespectrum. Then, applying the infinite loop functor $\Omega^{\infty}$, we obtain the following diagram of homotopy fiber sequences. 
%%%Taking the fibrant and using the lemma p.288 [MP].
\begin{center}
\begin{equation}\label{Diag:Bigdiagramspacelevel}
\begin{tikzpicture}[baseline = (current bounding box.center)]
 \node (AK-1) at (1,8) {$\Omega K_{a}$};
 \node (AQ-1) at (1,6) {$\Omega K_{a,\mathbb{Q}}$};
 \node (FAQ)  at (1,4) {$F_{a,\mathbb{Q}/\mathbb{Z}}$};
 \node (AK)   at (1,2) {$K_{a} $};
 \node (AQ)   at (1,0) {$K_{a,\mathbb{Q}}$};

 \node (TK1-1)at (4,8) {$\Omega K_{t}$};
 \node (TQ-1) at (4,6) {$\Omega K_{t,\mathbb{Q}}$};
 \node (FTQ)  at (4,4) {$F_{t,\mathbb{Q}/\mathbb{Z}}$};
 \node (TK1)  at (4,2) {$K_{t} $};
 \node (TQ)   at (4,0) {$K_{t,\mathbb{Q}}$};
 
 \node (TK2-1)at (7,8) {$\Omega K_{t}$};
 \node (TF-1) at (7,6) {$\Omega K_{t,\mathbb{F}^{\prime}}$};
 \node (FT)   at (7,4) {$F_{t,\mathbb{F}^{\prime}/\mathbb{Z}}$};
 \node (TK2)  at (7,2) {$K_{t}$};
 \node (TF)   at (7,0) {$K_{t,\mathbb{F}^{\prime}}$};

 \path [->, font=\scriptsize, >=angle 90]
 (AK-1) edge node [right]{$\Omega f$}(AQ-1)
 (AQ-1) edge (FAQ) 
 (FAQ)  edge (AK)
 (AK)   edge node [right]{$f$}(AQ)
 
 (TK1-1) edge node [right]{$\Omega g$}(TQ-1)
 (TQ-1)  edge (FTQ)
 (FTQ)   edge (TK1)
 (TK1)   edge node [right]{$g$}(TQ)
 
 (TK2-1)edge node [right]{$\Omega h$}(TF-1)
 (TF-1) edge (FT)
 (FT)   edge (TK2)
 (TK2)  edge node [right]{$h$}(TF)

 (AK-1) edge (TK1-1)
 (TK1-1)edge (TK2-1)
 (AQ-1) edge (TQ-1)
 (TQ-1) edge (TF-1) 
 (FAQ)  edge node [above]{$\operatorname{Su}$} node [below]{$\sim$}(FTQ)
 (FTQ)  edge (FT) 
 (AK)   edge (TK1)
 (TK1)  edge node [above]{$id$}(TK2)
 (AQ)   edge (TQ)
 (TQ)   edge (TF) ; 
  
\end{tikzpicture}
\end{equation}
\end{center} 
Now, since $\mathbf{K}_{t}\wedge \mathbf{M}\mathbb{Q}$ is rational, the following composition 
\[\mathbf{K}_{a}\wedge \mathbf{M}\mathbb{Z}\rightarrow \mathbf{K}_{t}\wedge \mathbf{M}\mathbb{Z} \rightarrow \mathbf{K}_{t} \wedge \mathbf{M}\mathbb{Q}\]
is determined by their induced homomorphisms (see \cite[Lemma $2.2.7$]{Wang1} or \cite[Theorem 5.8 and 7.11]{Ru}):
\begin{equation}\label{Eq:existenceofeh}
\pi_{\ast}(\mathbf{K}_{a}\wedge\mathbf{M}\mathbb{Z})\rightarrow \pi_{\ast}(\mathbf{K}_{t}\wedge\mathbf{M}\mathbb{Z})\rightarrow \pi_{\ast}(\mathbf{K}_{t}\wedge\mathbf{M}\mathbb{Q})=\pi_{\ast}(\mathbf{K}_{t})\otimes \mathbb{Q}.
\end{equation} 
Since, up to torsion groups, the algebraic $K$-groups of the real (resp. complex) numbers are divisible \cite[Theorem $4.9$]{Su}, \cite[VI.Theorem 1.6; Theorem 3.1]{We1} and there is no non-trivial homomorphism from a divisible group to a finitely generated abelian group, the composition \eqref{Eq:existenceofeh} must be zero. Thus, the composition $K_{a}\rightarrow K_{t}\rightarrow K_{t,\mathbb{Q}}$ is null-homotopic as an infinite loop map (or in $\operatorname{Ho}(\mathcal{P})$)\footnote{Another argument without using Suslin's result can be found in \cite[Proposition $2$]{We2}.}.

\begin{lemma}\phantomsection\label{Existanceofeh}
Liftings of the comparison map $K_{a}\rightarrow K_{t}$ with respect to the homotopy fiber sequence $F_{t,\mathbb{F/Z}}\rightarrow K_{t}\rightarrow K_{t,\mathbb{F}}$ exist, denoted by $e_{h}^{\mathbb{F}}$, and they fit into the commutative diagram below. 
\begin{center}
\begin{tikzpicture}
 \node (AK-1) at (1,8) {$\Omega K_{a}$};
 \node (AQ-1) at (1,6) {$\Omega K_{a,\mathbb{Q}}$};
 \node (FAQ)  at (1,4) {$F_{a,\mathbb{Q}/\mathbb{Z}}$};
 \node (AK)   at (1,2) {$K_{a} $};
 \node (AQ)   at (1,0) {$K_{a,\mathbb{Q}}$};

 \node (TK1-1)at (4,8) {$\Omega K_{t}$};
 \node (TQ-1) at (4,6) {$\Omega K_{t,\mathbb{Q}}$};
 \node (FTQ)  at (4,4) {$F_{t,\mathbb{Q}/\mathbb{Z}}$};
 \node (TK1)  at (4,2) {$K_{t} $};
 \node (TQ)   at (4,0) {$K_{t,\mathbb{Q}}$};
 
 \node (TK2-1)at (7,8) {$\Omega K_{t}$};
 \node (TF-1) at (7,6) {$\Omega K_{t,\mathbb{F}^{\prime}}$};
 \node (FT)   at (7,4) {$F_{t,\mathbb{F}^{\prime}/\mathbb{Z}}$};
 \node (TK2)  at (7,2) {$K_{t}$};
 \node (TF)   at (7,0) {$K_{t,\mathbb{F}^{\prime}}$};

 \path [->, font=\scriptsize, >=angle 90]
 (AK-1) edge node [right]{$\Omega f $}(AQ-1)
 (AQ-1) edge (FAQ) 
 (FAQ)  edge (AK)
 (AK)   edge node [right]{$f $}(AQ)
 
 (TK1-1) edge node [right]{$\Omega g $}(TQ-1)
 (TQ-1)  edge (FTQ)
 (FTQ)   edge (TK1)
 (TK1)   edge node [right]{$g $}(TQ)
 
 (TK2-1)edge node [right]{$\Omega h $}(TF-1)
 (TF-1) edge (FT)
 (FT)   edge (TK2)
 (TK2)  edge node [right]{$h $}(TF)

 (AK-1) edge (TK1-1)
 (TK1-1)edge (TK2-1)
 (AQ-1) edge (TQ-1)
 (TQ-1) edge (TF-1) 
 (FAQ)  edge node [above]{$\operatorname{Su}$}node[below]{$\sim$}(FTQ)
 (FTQ)  edge node [above]{$j$} (FT) 
 (AK)   edge (TK1)
 (TK1)  edge node [above]{$id$}(TK2)
 (AQ)   edge (TQ)
 (TQ)   edge (TF) ;

 \draw [dashed, ->] (AK) to [out=40,in=-120] node [yshift=.7em, xshift=1em]{\scriptsize $e_{h}^{\mathbb{F}^{\prime}}$}(FT);
 \draw [dashed, ->] (AK) to [out=45,in=-120] node [above]{\scriptsize $e_{h}^{\mathbb{Q}}$}(FTQ); 
 
\end{tikzpicture}
\end{center}
\end{lemma}
\begin{proof}
As shown in the discussion preceding the lemma, the composition $K_{a}\rightarrow K_{t}\rightarrow K_{t,\mathbb{F}}$ is null-homotopic, and choosing a null-homotopy gives us a lifting $e_{h}^{\mathbb{F}}$. To see the diagram is commutative, we note $F_{a,\mathbb{Q/Z}}$ has all its homotopy groups are torsion groups. That is because both $\operatorname{coker}(\pi_{n}(\Omega f ))$ and $\operatorname{ker}(\pi_{n}(f ))$ are torsion groups, for $n\geq 1$, and $\pi_{n}(F_{a,\mathbb{Q/Z}})$ fits into the short exact sequence 
\[0\rightarrow \operatorname{coker}(\pi_{n}(\Omega f)) \rightarrow \pi_{n}(F_{a,\mathbb{Q/Z}})\rightarrow \operatorname{ker}(\pi_{n}(f))\rightarrow 0.\] Therefore the rationalization of $F_{a,\mathbb{Q/Z}}$ is contractible and the abelian group $[F_{a,\mathbb{Q/Z}},\Omega K_{t,\mathbb{F}}]_{\operatorname{Ho}(\mathcal{P})}$ is trivial. This implies the commutativity of the following two triangles:   
\begin{center}
\begin{tikzpicture}
\node(Lu) at (0,2) {$F_{a,\mathbb{Q/Z}}$};
\node(Ll) at (0,0) {$K_{a}$}; 
\node(Ru) at (3,2) {$F_{t,\mathbb{Q/Z}}$};
 
\path[->, font=\scriptsize,>=angle 90]

(Lu) edge node [above]{$\operatorname{Su}$} (Ru)  
(Lu) edge (Ll);
 \draw [dashed, ->] (Ll) to [out=0,in=-90] node [yshift=.7em]{\scriptsize $e_{h}^{\mathbb{Q}}$}(Ru);

\node(Lu) at (5,2) {$F_{a,\mathbb{Q/Z}}$};
\node(Ll) at (5,0) {$K_{a}$}; 
\node(Ru) at (8,2) {$F_{t,\mathbb{F^{\prime}/Z}}$};
 
\path[->, font=\scriptsize,>=angle 90]

(Lu) edge node [above]{$j\circ \operatorname{Su}$} (Ru)  
(Lu) edge (Ll);
 \draw [dashed, ->] (Ll) to [out=0,in=-90] node [yshift=.7em]{\scriptsize $e_{h}^{\mathbb{F}^{\prime}}$}(Ru); 
\end{tikzpicture}
\end{center}
Thus, the lemma is proved. 
\end{proof}  

\begin{theorem}\phantomsection\label{LiftingsandTor}
Given a pointed topological space $X$ and liftings $e_{h}^{\mathbb{Q}}$ and $e_{h}^{\mathbb{F}^{\prime}}$ as above, then they induce the following isomorphisms: 
\begin{align*}
e_{h,\ast}^{\mathbb{Q}}\vert_{\operatorname{Tor}}&: \operatorname{Tor}[X,K_{a}]\xrightarrow{\sim} [X,F_{t,\mathbb{Q}/\mathbb{Z}}]\\
e_{h,\ast}^{\mathbb{F}^{\prime}}\vert_{\operatorname{Tor}}&:\operatorname{Tor}[X,K_{a}]\xrightarrow{\sim}  \operatorname{Tor}[X,F_{t,\mathbb{F^{\prime}}/\mathbb{Z}}]\simeq [X,F_{t,\mathbb{Q}/\mathbb{Z}}].
\end{align*}

\end{theorem}

\begin{proof}
\textbf{Step $1:$} We want to show $j_{\ast}$ induces an isomorphism  
\[\operatorname{coker}((\Omega g)_{\ast})\rightarrow \operatorname{Tor}((\operatorname{coker}(\Omega h)_{\ast}).\] 
To see this, we first recall some facts in homological algebra: Given an abelian group $A$, there is a short exact sequence
 \[0\rightarrow A_{T}\rightarrow A\rightarrow A/A_{T}\rightarrow 0,\]
where $A_{T}$ is the torsion subgroup of $A$. Now, by the right exactness of the tensor product, and the fact that $A_{T}\otimes \mathbb{F/Z}=0$, we see the homomorphism
\begin{equation}\label{Eq1inhomotopyliftings}
A\otimes \mathbb{F/Z}\rightarrow A/A_{T}\otimes \mathbb{F/Z}
\end{equation}
is an isomorphism. Moreover, since $A/A_{T}$ is torsion free and hence flat, there is another short exact sequence 
\[0\rightarrow A/A_{T}\otimes \mathbb{Q/Z}\rightarrow A/A_{T}\otimes \mathbb{F}^{\prime}/\mathbb{Z}\rightarrow A/A_{T}\otimes \mathbb{F}^{\prime}/\mathbb{Q}\rightarrow 0.\]
Because $A/A_{T}$ and $\mathbb{F}^{\prime}/\mathbb{Q}$ both are flat, the tensor product $A/A_{T}\otimes \mathbb{F}^{\prime}/\mathbb{Q}$ is flat and hence torsion free. Applying the left exactness of $\operatorname{Tor}$, we further obtain the following isomorphism 
\begin{equation}\label{Eq2inhomotopyliftings}
A/A_{T}\otimes \mathbb{Q/Z}\rightarrow \operatorname{Tor}(A/A_{T}\otimes \mathbb{F}^{\prime}/\mathbb{Z}).
\end{equation}  
Return to the theorem and let $A$ be the abelian group $[X,\Omega K_{t}]$. Then we see the isomorphisms \eqref{Eq1inhomotopyliftings} and \eqref{Eq2inhomotopyliftings} give the isomorphism claimed  
\[\operatorname{coker}((\Omega g )_{\ast})\cong A/A_{T}\otimes \mathbb{Q/Z}\xrightarrow{\sim} \operatorname{Tor}(A/A_{T}\otimes \mathbb{F}^{\prime}/\mathbb{Z})\cong \operatorname{Tor}(\operatorname{coker}((\Omega h)_{\ast}).\]

\textbf{Step $2$:} We claim that the homomorphism
\[[X,F_{t,\mathbb{Q/Z}}]\rightarrow \operatorname{Tor}[X,F_{t,\mathbb{F/Z}}]\]
is an isomorphism. This can be seen from the following two diagrams of exact sequences: The first one is obtained from the diagram \eqref{Diag:Bigdiagramspacelevel}:
\begin{center}
\begin{tikzpicture}
\node (01) at (0,1){$0$};
\node (00) at (0,0){$0$};
\node (11) at (2,1){$\operatorname{coker}((\Omega g )_{\ast})$};
\node (10) at (2,0){$\operatorname{coker}((\Omega h )_{\ast})$};
\node (21) at (5,1){$[X,F_{t,\mathbb{Q/Z}}]$};
\node (20) at (5,0){$[X,F_{t,\mathbb{F}^{\prime}/\mathbb{Z}}]$};
\node (31) at (8,1){$\operatorname{ker}(g_{ \ast})$};
\node (30) at (8,0){$\operatorname{ker}(h_{ \ast})$};
\node (41) at (10,1){$0$};
\node (40) at (10,0){$0,$};

 \path [->, font=\scriptsize, >=angle 90]

(00) edge (10)
(10) edge (20)
(20) edge (30)
(30) edge (40)

(01) edge (11)
(11) edge (21)
(21) edge (31)
(31) edge (41)
(11) edge (10)
(21) edge (20);

\draw [double equal sign distance] (31) to (30);  
\end{tikzpicture}
\end{center}
Applying the functor $\operatorname{Tor}$, we obtain the second one: 
\begin{center}
\begin{tikzpicture}
\node (01) at (0,1){$0$};
\node (00) at (0,0){$0$};
\node (11) at (2,1){$\operatorname{coker}((\Omega g )_{\ast})$};
\node (10) at (2,0){$\operatorname{Tor}(\operatorname{coker}((\Omega h )_{\ast})$};
\node (21) at (5.7,1){$[X,F_{t,\mathbb{Q/Z}}]$};
\node (20) at (5.7,0){$\operatorname{Tor}([X,F_{t,\mathbb{F}^{\prime}/\mathbb{Z}}])$};
\node (31) at (8.5,1){$\operatorname{ker}(g_{\ast})$};
\node (30) at (8.5,0){$\operatorname{ker}(h_{\ast})$};
\node (41) at (10,1){$0$};
\node (40) at (10,0){$0$};

 \path [->, font=\scriptsize, >=angle 90]

(00) edge (10)
(10) edge (20)
(20) edge node [above]{$l$}(30)
(30) edge (40)

(01) edge (11)
(11) edge (21)
(21) edge (31)
(31) edge (41)

(11) edge (10)
(21) edge (20);

\draw [double equal sign distance] (31) to (30);  
\end{tikzpicture}
\end{center}
Note that the functor $\operatorname{Tor}$ does not always preserve short exact sequences, but, in this case, it does. The only thing to check is the surjectivity of $l$, yet it follows quickly from the surjectivity of the homomorphism 
\[[X,F_{t,\mathbb{Q/Z}}]\rightarrow \operatorname{Tor}([X,F_{t,\mathbb{F/Z}}])\rightarrow \operatorname{ker}(h_{\ast})=\operatorname{ker}(g_{\ast}).\]  
By the short five lemma, we see the homomorphism
\begin{equation}\label{Eq:IsoonTor}
[X,F_{t,\mathbb{Q/Z}}]\rightarrow \operatorname{Tor}[X,F_{t,\mathbb{F/Z}}]
\end{equation}
is indeed an isomorphism.

\textbf{Step $3$:} Consider the commutative diagram below:

\begin{center}
\begin{tikzpicture}
\node (00) at (0,0){$[X,K_{a}]$};
\node (03) at (0,3){$[X,F_{a,\mathbb{Q/Z}}]$};
\node (115) at (2,1.5){$\operatorname{Tor}[X,K_{a}]$};
\node (33) at (7,3){$[X,F_{t,\mathbb{Q/Z}}]$};

 \path [->, font=\scriptsize, >=angle 90]

(03) edge (00)
(03) edge node [above]{$t$} (115)
(115) edge (00)
(03) edge node [above]{$\operatorname{Su}_{\ast}$} node [below]{$\sim $}(33);

\draw [->](115) to [out=0,in=-110] node [above]{\scriptsize $e_{h,\ast}^{\mathbb{Q}}\vert_{\operatorname{Tor}}$}(33);
\draw [->](00) to [out=0,in=-90] node [above]{\scriptsize $e_{h,\ast}^{\mathbb{Q}}$} (33);
\end{tikzpicture}
\end{center}
and note the homomorphism $[X,F_{a,\mathbb{Q/Z}}]\rightarrow [X,K_{a}]$ factors through the homomorphism $t$. Since $\operatorname{Su}_{\ast}$ is an isomorphism by Suslin's theorem (see \cite[Theorem $2.1.8$]{Wang1}), we know $t$ is injective. On the other hand, $t$ is surjective by its definition, and thus $t$ is an isomorphism. Now, in view of the diagram above, we see the homomorphism $e_{h,\ast}^{\mathbb{Q}}\vert_{\operatorname{Tor}}$ has to be an isomorphism as well.

Combining with the isomorphism \eqref{Eq:IsoonTor}, the second assertion can be deduced from the following commutative diagram:
\begin{center}
\begin{tikzpicture}
\node (00) at (0,0){$[X,K_{a}]$};
\node (03) at (0,3){$[X,F_{a,\mathbb{Q/Z}}]$};
\node (115) at (2,1.5){$\operatorname{Tor}[X,K_{a}]$};
\node (33) at (3,3){$[X,F_{t,\mathbb{Q/Z}}]$};
\node (63) at (6,3){$\operatorname{Tor}([X,F_{t,\mathbb{F}^{\prime}/\mathbb{Z}}])$};
\node (93) at (9,3){$[X,F_{\mathbb{F}^{\prime}/\mathbb{Z}}]$};

 \path [->, font=\scriptsize, >=angle 90]

(03) edge (00)
(03) edge node [above]{$t$} (115)
(115) edge  (00)
 
 ;
\draw (33) to node [above]{\scriptsize $j_{\ast}$} node[below]{$\sim$}(63);
\draw (03) to node [above]{\scriptsize $\operatorname{Su}_{\ast}$}node [below]{$\sim$}(33);
\draw [right hook->] (63) to (93);

\draw [->](115) to [out=0,in=-100] node [yshift=1em]{\scriptsize $e_{h,\ast}^{\mathbb{F}^{\prime}}\vert_{\operatorname{Tor}}$}(63);
\draw [->](00) to [out=0,in=-90] node [above]{\scriptsize $e_{h,\ast}^{\mathbb{F}^{\prime}}$} (93);
\end{tikzpicture}
\end{center} 
\end{proof}

In fact, we can see from the diagrams above the homomorphisms $e_{h,\ast}^{\mathbb{F}^{\prime}}\vert_{\operatorname{Tor}}$ and $e_{h,\ast}^{\mathbb{Q}}\vert_{\operatorname{Tor}}$ are identical to the compositions $j_{\ast}\circ \operatorname{Su}_{\ast}\circ t^{-1}$ and $\operatorname{Su}_{\ast}\circ t^{-1}$, respectively. This means any lifting induces the same isomorphism on the torsion subgroup $\operatorname{Tor}[X,K_{a}]$. On the other hand, when $X=L$, a finite $CW$-complex, the $e$-invariant satisfies the same commutative diagram \cite[Lemma $4.4.2$]{Wang1}: 
\begin{center}
\begin{tikzpicture}
\node(Lu) at (0,2) {$[L,F_{a,\mathbb{Q/Z}}]$};
\node(Ll) at (0,0) {$[L,K_{a}\mathbb{C}]$}; 
\node(Ru) at (3,2) {$[L,F_{t,\mathbb{Q/Z}}]$};
\node(RRu) at (6,2) {$[L,F_{t,\mathbb{C/Z}}]$};

\path[->, font=\scriptsize,>=angle 90]

(Lu) edge node [above]{$\operatorname{Su}_{\ast}$}node [below]{$\sim$}(Ru)  
(Lu) edge (Ll)  
(Ru) edge (RRu);
\draw[->] (Ll) to [out=0,in=-100] node[above]{\scriptsize $e_{\ast}$}(RRu);
\end{tikzpicture}
\end{center}
Hence, we have shown the following theorem: 
\begin{theorem}\phantomsection\label{sameisoLiftings}
Every lifting $e_{h}^{\mathbb{F}}:K_{a}\rightarrow F_{t,\mathbb{F/Z}}$ of the comparison map $K_{a}\rightarrow K_{t}$ induces the same isomorphism 
\[e_{h,\ast}^{\mathbb{F}}\vert_{\operatorname{Tor}}: \operatorname{Tor}[X,K_{a}]\xrightarrow{\sim} \operatorname{Tor}[X,F_{t,\mathbb{F}/\mathbb{Z}}].\]

When $X=L$, a finite $CW$-complex, and in the case of the complex numbers, $e_{h,\ast}^{\mathbb{F}}$ restricts to the $e$-invariant on the torsion subgroup of $[L,K_{a}\mathbb{C}]$.
\end{theorem}

\section{Infinitely many different liftings}
In the last section we have seen the existence of liftings of the comparison map  
\[K_{a}\rightarrow K_{t}.\]
We shall see in this section, in the case of the complex numbers, there are infinitely many different liftings of the comparison map $K_{a}\mathbb{C}\rightarrow K_{t}$ in the category $\operatorname{Ho}(\mathcal{P})$.

Recall that the number of different liftings is measured by the size of the subgroup 
\begin{equation}\label{Eq:Countingliftings}
\operatorname{Im}([K_{a}\mathbb{C},\Omega K_{t,\mathbb{F}}]_{\operatorname{Ho}(\mathcal{P})})\subset
[K_{a}\mathbb{C},F_{t,\mathbb{F/Z}}]_{\operatorname{Ho}(\mathcal{P})}.
\end{equation} 
This results from the following long exact sequence
\[...\rightarrow [K_{a}\mathbb{C},\Omega K_{t,\mathbb{F}}]_{\operatorname{Ho}(\mathcal{P})}\rightarrow [K_{a}\mathbb{C},F_{t,\mathbb{F/Z}}]_{\operatorname{Ho}(\mathcal{P})}\rightarrow [K_{a}\mathbb{C}.K_{t}]_{\operatorname{Ho}(\mathcal{P})}\rightarrow....\]
Since $\Omega K_{t,\mathbb{F}}$ is rational, by \cite[Lemma $2.2.7$]{Wang1}, we have the isomorphism of abelian groups
\begin{equation}\label{Eq:Infinitemanyliftings}
[K_{a}\mathbb{C},\Omega K_{t,\mathbb{F}}]_{\operatorname{Ho}(\mathcal{P})}\xrightarrow{\sim} \operatorname{Hom}^{0}(\pi_{\ast}(K_{a}\mathbb{C})\otimes\mathbb{Q},\pi_{\ast}(\Omega K_{t,\mathbb{F}})),
\end{equation}  
where $\operatorname{Hom}^{0}(A_{\ast},B_{\ast})$ is the abelian group of homogeneous homomorphisms of degree $0$ between graded abelian groups $A_{\ast}$ and $B_{\ast}$. Now it is known that the abelian group 
\[
  \pi_{\ast}(\Omega K_{t,\mathbb{F}})=
              \begin{cases}
               \mathbb{F} & \ast=odd\\
               0& \ast=even,\\ 
               \end{cases}
\] 
while, to the author's knowledge, the precise size of the abelian group $\pi_{\ast}(K_{a}\mathbb{C})\otimes \mathbb{Q}$ is not determined yet. Nevertheless, according to \cite[Sec.4-5]{Jah}, we have $\pi_{\ast}(K_{a}\mathbb{C})\otimes\mathbb{Q}$ is a non-trivial $\mathbb{Q}$-vector space, when $\ast$ is odd. In fact, Jahren constructs a homomorphism from $\pi_{\ast}(K_{a}\mathbb{C})\rightarrow \mathbb{R}$, for $\ast$ is odd and proves that this homomorphism reduces to the Borel classes after precomposing the homomorphisms induced by the conjugate embeddings of a number field in $\mathbb{C}$ and tensoring $\mathbb{R}$. Hence, when $\ast$ is odd, $\pi_{\ast}(K_{a}\mathbb{C})\otimes \mathbb{Q}$ cannot be trivial.

In view of \eqref{Eq:Countingliftings} and \eqref{Eq:Infinitemanyliftings}, we know if one can construct infinitely many different homomorphisms
\[\pi_{\ast}(K_{a}\mathbb{C})\rightarrow \pi_{\ast}(\Omega K_{t,\mathbb{F}})\]
such that, after composing with the homomorphism
\[\pi_{\ast}(\Omega K_{t,\mathbb{F}})\rightarrow \pi_{\ast}(F_{t,\mathbb{F/Z}}),\]
they remain different, then we obtain infinitely many different liftings. We provide one construction here: Let $\ast$ be an odd number and pick up a non-trivial element $x\in \pi_{\ast}(K_{a}\mathbb{C})\otimes \mathbb{Q}$. Assign to it the number $\frac{1}{n}\in \mathbb{F}$ with $n\in\mathbb{N}\setminus \{1\}$. Then extend this assignment to a homomorphism 
\[\pi_{\ast}(K_{a}\mathbb{C})\otimes \mathbb{Q}\rightarrow \mathbb{F}=\pi_{\ast}(\Omega F_{t,\mathbb{F}}).\]
It is not to difficult to find an extension of this assignment. For instance, one can choose an inner product on $\pi_{\ast}(K_{a}\mathbb{C})\otimes\mathbb{Q}$ and let $<x>^{\perp}$ go to zero. Therefore, we have shown that the subgroup 
\[\operatorname{Im}([K_{a}\mathbb{C},\Omega K_{t,\mathbb{F}}]_{\operatorname{Ho}(\mathcal{P})})\subset [K_{a}\mathbb{C},F_{t,\mathbb{F/Z}}]_{\operatorname{Ho}(\mathcal{P})}\] 
contains at least countably infinitely many different elements.
\begin{theorem}\phantomsection\label{Manyliftings}
There are infinitely many different liftings $K_{a}\mathbb{C}\rightarrow F_{t,\mathbb{F/Z}}$ of the comparison map $K_{a}\mathbb{C}\rightarrow K_{t}$ such that the following diagram commutes:

\begin{center}
\begin{tikzpicture}
\node(Lu) at (0,2) {$F_{a,\mathbb{Q/Z}}$};
\node(Ll) at (0,0) {$K_{a}\mathbb{C}$}; 
\node(Ru) at (2,2) {$F_{t,\mathbb{F/Z}}$};
\node(Rl) at (2,0) {$K_{t}$};

\path[->, font=\scriptsize,>=angle 90]

(Lu) edge (Ru)  
(Lu) edge (Ll)
(Ll) edge (Rl) 
(Ru) edge (Rl);

\draw[dashed,->] (Ll) to (Ru);

\end{tikzpicture}

\end{center}
\end{theorem}
\begin{proof}
This results from the discussion preceding the theorem and the fact that the rationalization of $F_{a,\mathbb{Q/Z}}$ is contractible. 
\end{proof}

\begin{remark}
The same method does not work in the case of the real numbers. In fact, by the theorem of Jahren \cite{Jah}, we can only conclude $\pi_{\ast}(K_{a}\mathbb{R})\otimes \mathbb{Q}$ is non-trivial when $\ast=4k-3$. On the other hand, we have $\pi_{\ast}(\Omega K_{t,\mathbb{F}})=\mathbb{F}$ when $\ast=4k-1$ and $0$ otherwise, where $k\in\mathbb{N}$.
\end{remark}

\section{The maps $e$, $\operatorname{ch}^{rel}$, $e_{h}^{\natural}$ and $t^{\natural}_{h}$}
In this section, we construct the maps $e_{h}^{\natural}$ and $t_{h}^{\natural}$ and study their relation with the maps $e$ and $\operatorname{ch}^{rel}$ in \cite[Sections $4.1$ and $4.3$]{Wang1}. We begin with some structure theorems for $K_{a}$ and $K^{rel}$.  

\begin{corollary}\label{StrKa}
There exists a homotopy equivalence of infinite loop spaces:
\[K_{a}\xrightarrow{\sim} K_{a,\mathbb{Q}}\times F_{t,\mathbb{Q/Z}}.\] 
\end{corollary}
\begin{proof}
By Lemma \ref{Existanceofeh}, there exists a lifting $e_{h}^{\mathbb{Q}}: K_{a}\rightarrow F_{t,\mathbb{Q/Z}}$ of the comparison map $K_{a}\rightarrow K_{t}$. Combining with the rationalization $u_{\mathbb{Q}}:K_{a}\rightarrow K_{a,\mathbb{Q}}$, we obtain a homotopy equivalence of infinite loop spaces
\[K_{a}\xrightarrow{(u_{\mathbb{Q}},e_{h}^{\mathbb{Q}})} K_{a,\mathbb{Q}}\times F_{t,\mathbb{Q/Z}}.\]
\end{proof}

\begin{corollary}\phantomsection\label{Torsionpartarethesame}
Given a lifting $e_{h}^{\mathbb{F}}:K_{a}\rightarrow F_{t,\mathbb{F/Z}}$, the composition 
\[F_{t,\mathbb{Q/Z}}\xrightarrow{i_{2}} K_{a,\mathbb{Q}}\times F_{t,\mathbb{Q/Z}}\simeq K_{a}\xrightarrow{e_{h}^{\mathbb{F}}} F_{t,\mathbb{F/Z}}\]
is homotopic, as an infinite loop map, to the canonical map
\[j:F_{t,\mathbb{Q/Z}}\rightarrow F_{t,\mathbb{F/Z}}\]
induced by the inclusion $\mathbb{Q/Z}\hookrightarrow \mathbb{C/Z},$ where $i_{2}$ is the inclusion into the second component and $K_{a}\simeq K_{a,\mathbb{Q}}\times F_{t,\mathbb{Q/Z}}$ is the homotopy equivalence given by a lifting $e^{\mathbb{Q}}_{h}$. 

In other words, what determines a lifting $e_{h}^{\mathbb{F}}$ is its restriction to the divisible part $K_{a,\mathbb{Q}}$.

\end{corollary} 
\begin{proof}
This follows from the commutative diagram below
\begin{center}
\begin{equation}\label{Diag:OnlyKaQmatters1} 
\begin{tikzpicture}[baseline=(current bounding box.center)]
 \node (Lu) at (0,2) {$F_{t,\mathbb{Q/Z}}$};
 \node (Ll)  at (0,0) {$K_{a,\mathbb{Q}}\times F_{t,\mathbb{Q/Z}}$};
 \node (Mu)   at (4,2) {$F_{a,\mathbb{Q/Z}}$};
 \node (Ml)   at (4,0) {$K_{a}$};
 \node (Rl)   at (8,0) {$F_{t,\mathbb{F/Z}}$};
 \node (hty) at (3.8,-.6) {$\wr\wr$};

 \path [->, font=\scriptsize, >=angle 90]
  (Lu)  edge node [right]{$i_{2}$} (Ll)
  (Mu)  edge node [above]{$Su\hspace{1em}$} node[below]{$\sim$} (Lu) 
  (Ml)  edge node [above]{$(u_{\mathbb{Q}},e_{h}^{\mathbb{Q}})$} node [below]{$\sim$}(Ll) 
  (Mu)  edge (Ml)
  (Ml)  edge node [above]{$e_{h}^{\mathbb{F}}$}(Rl);
  \draw [->] (Lu) to [out=40,in=130] node [below]{\scriptsize $j$}(Rl);
  \draw [->] (Ll) to [out=-30,in=-159] (Rl); 
\end{tikzpicture} 
\end{equation}
\end{center}
\end{proof}

Recall that $K^{rel}$ is the homotopy fiber of the comparison map $K_{a}\rightarrow K_{t}$. 

\begin{corollary}\label{StrKr}
There exists a homotopy equivalence of infinite loop spaces:
\[K^{rel}\xrightarrow{\sim} K_{a,\mathbb{Q}}\times  \Omega  K_{t,\mathbb{Q}}.\]
\end{corollary}
\begin{proof}
Choose a lifting $e_{h}^{\mathbb{Q}}$ and hence a homotopy equivalence $K_{a}\xrightarrow{(u_{\mathbb{Q}},e_{h}^{\mathbb{Q}})}K_{a,\mathbb{Q}}\times F_{t,\mathbb{Q/Z}}$, and consider the commutative diagram below:
\begin{center}
\begin{equation}\label{Diag:Relationbewteeneandh}
\begin{tikzpicture}[baseline=(current bounding box.center)]
 \node (OKt) at (0,6) {$\Omega K_{t}$};
 \node (Kr)  at (0,4) {$K^{rel}$};
 \node (Ka)   at (0,2) {$K_{a} $};
 \node (Kt)   at (0,0) {$K_{t}$};

 \node (OKt1) at (4,6) {$\Omega K_{t}$};
 \node (Pro2)  at (4,4) {$K_{a,\mathbb{Q}}\times \Omega K_{t,\mathbb{Q}}$};
 \node (Pro1)  at (4,2) {$K_{a,\mathbb{Q}}\times F_{t,\mathbb{Q/Z}} $};
 \node (Kt1)   at (4,0) {$K_{t}$};

 \path [->, font=\scriptsize, >=angle 90]
  (Kr)  edge node [above]{$(u_{\mathbb{Q}}\circ \pi,t^{\mathbb{Q}}_{h})$} (Pro2)
  (Ka)  edge node [above]{$(u_{\mathbb{Q}},e_{h}^{\mathbb{Q}})$} (Pro1) 
  
  (OKt)  edge (Kr)
  (Kr)   edge node [right]{$\pi$}(Ka)
  (Ka)   edge (Kt)
  (OKt1) edge (Pro2) 
  (Pro2) edge  (Pro1)
  (Pro1) edge node [right]{$p$}(Kt1); 
  \draw [double equal sign distance](OKt) to (OKt1);
  \draw [double equal sign distance] (Kt) to (Kt1); 
\end{tikzpicture} 
\end{equation}
\end{center}
where $p$ is the composition 
\[K_{a,\mathbb{Q}}\times F_{t,\mathbb{Q/Z}}\xrightarrow{\pi_{2}} F_{t,\mathbb{Q/Z}}\rightarrow K_{t},\]
$\pi_{2}$ is the projection onto the second component and $t^{\mathbb{Q}}_{h}$ is an infinite loop map induced by $e_{h}^{\mathbb{Q}}$ and a filler (homotopy) of the following triangle:
 
\begin{center}
\begin{equation}\label{Eq1:3.3}
\begin{tikzpicture}[baseline= (current bounding box.center)]
\node (Ll) at (0,2) {$K_{a}$};
\node (Ru) at (2,2) {$F_{t,\mathbb{Q/Z}}$};
\node (Rl) at (2,0) {$K_{t}$};

\draw [->](Ll) to node [above]{\scriptsize $e_{h}^{\mathbb{Q}}$}(Ru); 
\draw [->](Ru) to (Rl);
\draw [->](Ll) to (Rl);
\end{tikzpicture}
\end{equation}
\end{center}

Since $(u_{\mathbb{Q}},e_{h}^{\mathbb{Q}}):K_{a}\rightarrow K_{a,\mathbb{Q}}\times F_{t,\mathbb{Q/Z}}$ is a homotopy equivalence, in view of the commutative diagram \eqref{Diag:Relationbewteeneandh}, we see $(u_{\mathbb{Q}}\circ \pi,t_{h}^{\mathbb{Q}})$ is also a homotopy equivalence.
 
\end{proof} 

The next lemma shows what determines $t_{h}^{\mathbb{F}}$ is its restriction to $K_{a,\mathbb{Q}}$.
\begin{lemma}  
Let $t_{h}^{\mathbb{F}}:K^{rel}\rightarrow \Omega K_{t,\mathbb{F}}$ be a lifting of the composition $e_{h}^{\mathbb{Q}}\circ \pi:K^{rel}\rightarrow K_{a}\rightarrow F_{t,\mathbb{F/Z}}$ with respect to
the fiber sequence 
\[\Omega K_{t,\mathbb{F}}\rightarrow F_{t,\mathbb{F/Z}}\rightarrow K_{t}.\] 
Then in $\operatorname{Ho}(\mathcal{P})$ the composition 
\[l:\Omega K_{t,\mathbb{Q}}\xrightarrow{i_{2}} K_{a,\mathbb{Q}}\times \Omega K_{t,\mathbb{Q}}\rightarrow \Omega K_{t,\mathbb{F}}\]
is homotopic to the canonical map $\Omega K_{t,\mathbb{Q}}\xrightarrow{j^{\prime}} \Omega K_{t,\mathbb{F}}$ induced from the inclusion $\mathbb{Q}\hookrightarrow \mathbb{F}$.    
\end{lemma}
\begin{proof}
This amounts to show the commutative diagram below:
\begin{center}
\begin{tikzpicture}
 \node (Lu) at (0,2) {$\Omega K_{t,\mathbb{Q}}$};
 \node (Ll) at (0,0) {$K_{a,\mathbb{Q}}\times \Omega K_{t,\mathbb{Q}}$};
 \node (Ml)   at (4,0) {$K^{rel}$};
 \node (Rl)   at (8,0) {$\Omega K_{t,\mathbb{F}}$};
 \node (hty) at (3.8,-.6) {$\wr\wr$};

 \path [->, font=\scriptsize, >=angle 90]
  (Lu)  edge node [right]{$i_{2}$} (Ll)
  (Ml)  edge node [above]{$(u_{\mathbb{Q}}\circ \pi,t_{h}^{\mathbb{Q}})$} node [below]{$\sim$}(Ll) 
  (Ml)  edge node [above]{$t_{h}^{\mathbb{F}}$}(Rl);
  \draw [->] (Lu) to [out=0,in=150] node [below]{\scriptsize $j^{\prime}$}(Rl);
  \draw [->] (Ll) to [out=-30,in=-159] (Rl); 
\end{tikzpicture} 
\end{center}

Combining with the diagram \eqref{Diag:OnlyKaQmatters1}, we see that $l$ and $j$ are homotopic after composing with the map $\Omega K_{t,\mathbb{F}}\rightarrow F_{t,\mathbb{F/Z}}$. It means $l$ and $j$ differ by a map $\Omega K_{t,\mathbb{Q}}\rightarrow \Omega K_{t,\mathbb{F}}$ that factors through $\Omega K_{t}$. However, given any map $\Omega K_{t,\mathbb{Q}}\rightarrow \Omega K_{t}$, the induced homomorphism $\pi_{\ast}(\Omega K_{t,\mathbb{Q}})\rightarrow \pi_{\ast}(\Omega K_{t})$ is always trivial and hence the composition $\pi_{\ast}(\Omega K_{t,\mathbb{Q}})\rightarrow \pi_{\ast}(\Omega K_{t})\rightarrow \pi_{\ast}(\Omega K_{t,\mathbb{F}})$ is also trivial. Since $\Omega K_{t,\mathbb{F}}$ is rational, we see any map $\Omega K_{t,\mathbb{Q}}\rightarrow \Omega K_{t,\mathbb{F}}$ that factors through $\Omega K_{t}$ is null-homotopic \cite[Lemma $2.2.7$]{Wang1}. Therefore, $l$ and $j$ have to be homotopic in $\operatorname{Ho}(\mathcal{P})$.
\end{proof}

A priori, the map $t_{h}^{\mathbb{F}}$ depends on the choice of fillers of the diagram \eqref{Eq1:3.3}, the following shows, in effect, every filler induces the same $t_{h}^{\mathbb{F}}$ in $\operatorname{Ho}(\mathcal{P})$.
  
\begin{corollary}\label{edeterminest}
Given a lifting $e^{\mathbb{F}}_{h}$, there is a unique map 
\[t_{h}^{\mathbb{F}}:K^{rel}\rightarrow \Omega K_{t,\mathbb{F}}\] 
making the following diagram commute:
\begin{equation}\label{tecommutativediag}
\begin{tikzpicture}[baseline=(current  bounding  box.center)] 
 \node (OKt) at (0,6) {$\Omega K_{t}$};
 \node (Kr)  at (0,4) {$K^{rel}$};
 \node (Ka)   at (0,2) {$K_{a}$};
 \node (Kt)   at (0,0) {$K_{t}$};

 \node (OKt1) at (4,6) {$\Omega K_{t}$};
 \node (Pro2)  at (4,4) {$\Omega K_{t,\mathbb{F}}$};
 \node (Pro1)  at (4,2) {$F_{t,\mathbb{F/Z}}$};
 \node (Kt1)   at (4,0) {$K_{t}$};

 \path [->, font=\scriptsize, >=angle 90]
  (Kr)  edge node [above]{$t_{h}^{\mathbb{F}}$} (Pro2)
  (Ka)  edge node [above]{$e_{h}^{\mathbb{F}}$} (Pro1) 
  
  (OKt)  edge node [right]{$i$}(Kr)
  (Kr)   edge node [right]{$\pi$}(Ka)
  (Ka)   edge (Kt)
  (OKt1) edge (Pro2) 
  (Pro2) edge (Pro1)
  (Pro1) edge (Kt1); 
  \draw [double equal sign distance](OKt) to (OKt1);
  \draw [double equal sign distance] (Kt) to (Kt1); 
\end{tikzpicture} 
\end{equation}
\end{corollary}
\begin{proof}
Suppose there is another map $t_{h}^{\mathbb{F},\prime}$ which also fits into the commutative diagram \eqref{tecommutativediag}, then the difference between $t_{h}^{\mathbb{F}}$ and $t_{h}^{\mathbb{F},\prime}$ is measured by an element in the image  
\begin{equation}\label{Eq:Unieqnessofth}
\operatorname{Im}([K^{rel},\Omega K_{t}]_{\operatorname{Ho}(\mathcal{P})})\subset [K^{rel},F_{t,\mathbb{F/Z}}]_{\operatorname{Ho}(\mathcal{P})}.
\end{equation}
Since any map $K^{rel}\rightarrow \Omega K_{t}$ 
induces the trivial homomorphism between homotopy groups, the subgroup \eqref{Eq:Unieqnessofth} should be trivial, and hence, the maps $t_{h}^{\mathbb{F}}$ and $t_{h}^{\mathbb{F},\prime}$ must be homotopic in $\operatorname{Ho}(\mathcal{P})$. 
\end{proof} 

Now if we fix the map $t_{h}^{\mathbb{F}}$ in the diagram \eqref{tecommutativediag} instead, we get a slightly weaker result. Recall first that, given $\mathbf{E}$ a $CW$-prespectrum and $\mathbf{F}$ an $\Omega$-prespectrum, by a phantom map $\bm{f}:\mathbf{E}\rightarrow \mathbf{F}$ in $\mathcal{P}$, we understand its restriction to any finite $CW$-subprespectrum is null-homotopic.
\begin{lemma}\phantomsection\label{thdetermineseh}
The map $t_{h}^{\mathbb{F}}$ determines the lifting $e_{h}^{\mathbb{F}}$ up to phantom maps. Namely, if there is another lifting $e_{h}^{\mathbb{F},\prime}$ such that the pair $(t_{h}^{\mathbb{F}},e_{h}^{\mathbb{F},\prime})$ also satisfies the commutative diagram \eqref{tecommutativediag}, then $e_{h}^{\mathbb{F},\prime}$ and $e_{h}^{\mathbb{F}}$ differ by a phantom map.
\end{lemma}
\begin{proof}
Firstly, recall that $e_{h}^{\mathbb{Q}}$ induces an identification \[K_{a}\simeq K_{a,\mathbb{Q}}\times F_{t,\mathbb{Q/Z}},\]
and, by Corollary \ref{Torsionpartarethesame}, we know $e_{h}^{\mathbb{F}}$ and $e_{h}^{\mathbb{F},\prime}$ restrict to homotopic maps on $F_{t,\mathbb{Q/Z}}$ in $\operatorname{Ho}(\mathcal{P})$. 

Secondly, via the Serre class theory \cite[Proposition 4.23, 4.25]{Ru} (see \cite[A.2]{Wang1} for the relation between $\mathcal{P}$ and $\mathcal{A}$), one can deduce 
\[[\mathbf{F},\mathbf{K}_{t}]_{\operatorname{Ho}(\mathcal{P})}\] is a finitely generated abelian group, for every finite $CW$-prespectrum $\mathbf{F}$. On the other hand, since $\mathbf{K}_{a}\wedge \mathbf{M}\mathbb{Q}$ is rational, the abelian group 
\[[\mathbf{E},\mathbf{K}_{a}\wedge \mathbf{M}\mathbb{Q}]_{\operatorname{Ho}(\mathcal{P})}\]  
is always divisible, for any $CW$-prespectrum $\mathbf{E}$.
Hence the homomorphism 
\[ [\mathbf{F},\mathbf{K}_{a}\wedge \mathbf{M}\mathbb{Q}]_{\operatorname{Ho}(\mathcal{P})}\rightarrow [\mathbf{F},\mathbf{K}_{t}]_{\operatorname{Ho}(\mathcal{P})}\]
is trivial, for any finite $CW$-prespectrum $\mathbf{F}$. In particular, this implies all divisible elements of $[\mathbf{F},\mathbf{K}_{a}]_{\operatorname{Ho}(\mathcal{P})}$ are in the image of the homomorphism 
\[[\mathbf{F},\mathbf{K}^{rel}]_{\operatorname{Ho}(\mathcal{P})}\rightarrow [\mathbf{F},\mathbf{K}_{a}]_{\operatorname{Ho}(\mathcal{P})}.\] Using the commutative diagram \eqref{tecommutativediag} again, we conclude
\[e^{\mathbb{F}}_{h,\ast}=e^{\mathbb{F},\prime}_{h,\ast}:[\mathbf{F},\mathbf{K}_{a}]_{\operatorname{Ho}(\mathcal{P})}\rightarrow [\mathbf{F},\mathbf{Fib}(\mathbf{ch})]_{\operatorname{Ho}(\mathcal{P})},\] 
for every finite $CW$-prespectrum $\mathbf{F}$, where $\mathbf{Fib}(\mathbf{ch})$ is the homotopy fiber of 
\[\mathbf{K}_{t}\rightarrow \mathbf{K}_{t}\wedge \mathbf{M}\mathbb{F}.\]
Thus, we have proved the lemma.
\end{proof} 

We now compare $e_{h}^{\mathbb{F}}$ and $t_{h}^{\mathbb{F}}$ with the maps $e$ and $\operatorname{ch}^{rel}$ defined in \cite[Sections $4.1$ and $4.3$]{Wang1}. Recall the following lemma from \cite{Wang1}:
\begin{lemma}[Theorem $4.2.3$ in \cite{Wang1}]
The following diagram commutes up to weak homotopy
\begin{center}
\begin{equation} 
\begin{tikzpicture}[baseline=(current  bounding  box.center)] 
 \node (OKt) at (0,6) {$\Omega K_{t}$};
 \node (Kr)  at (0,4) {$K^{rel}\mathbb{C}$};
 \node (Ka)   at (0,2) {$K_{a}\mathbb{C}$};
 \node (Kt)   at (0,0) {$K_{t}$};

 \node (OKt1) at (4,6) {$\Omega K_{t}$};
 \node (Pro2)  at (4,4) {$\Omega K_{t,\mathbb{C}}$};
 \node (Pro1)  at (4,2) {$F_{t,\mathbb{C/Z}}$};
 \node (Kt1)   at (4,0) {$K_{t}$};

 \path [->, font=\scriptsize, >=angle 90]
  (Kr)  edge node [above]{$\operatorname{ch}^{rel}$} (Pro2)
  (Ka)  edge node [above]{$e$} (Pro1) 
  
  (OKt)  edge node [right]{$i$}(Kr)
  (Kr)   edge node [right]{$\pi$}(Ka)
  (Ka)   edge (Kt)
  (OKt1) edge node [right]{$\operatorname{ch}$}(Pro2) 
  (Pro2) edge (Pro1)
  (Pro1) edge (Kt1); 
  \draw [double equal sign distance](OKt) to (OKt1);
  \draw [double equal sign distance] (Kt) to (Kt1); 
\end{tikzpicture} 
\end{equation}
\end{center}
\end{lemma}
 
\begin{theorem}\phantomsection\label{enatural}
There is a lifting $e^{\natural}_{h}$ unique up to phantom maps such that  
\begin{align}
e^{\natural}_{h,\ast}=&e_{\ast}:\pi_{\ast}(K_{a}\mathbb{C})\rightarrow \pi_{\ast}(F_{t,\mathbb{C/Z}});\nonumber\\
e^{\natural}_{h,\ast}\vert_{\operatorname{Tor}}=&e_{\ast}\vert_{\operatorname{Tor}}:\operatorname{Tor}[L,K_{a}\mathbb{C}]\rightarrow [L,F_{t,\mathbb{C/Z}}];\label{Eq:Torenande}\\
\operatorname{ch}_{\otimes\mathbb{R}}\circ \operatorname{Im}\circ e^{\natural}_{h}=&\operatorname{Bo}\in \operatorname{Ho}(\mathcal{P}).\label{Eq:Boanden}
\end{align} 
%where $\operatorname{Bo}$ is the Borel regulator given in \cite[Theorem $3.1$]{JW} (see also \cite[]{Wang1}), $\operatorname{ch}$ is the Chern character $\Omega K_{t,\mathbb{R}}\rightarrow H^{odd}\mathbb{R}$ and $J:\Omega K_{t,\mathbb{C/Z}}\rightarrow \Omega K_{t,\mathbb{R}}$ is the map induced by taking the imaginary part of complex numbers \cite[]{Wang1}.   
\end{theorem}
\begin{proof}
Since $K_{t,\mathbb{Q}}$ is rational, there is a unique infinite loop map $t^{\natural}_{h}$ such that
\[t^{\natural}_{h,\ast}=\operatorname{ch}^{rel}_{\ast}:\pi_{\ast}(K^{rel}\mathbb{C})\rightarrow \pi_{\ast}(\Omega K_{t,\mathbb{C}}).\]
As the homotopy cofiber and fiber sequences are isomorphic in the stable homotopy category, choosing a filler of the following triangle
\begin{center} 
\begin{tikzpicture}
\node (Ll) at (0,0) {$K^{rel}\mathbb{C}$};
\node (Rl) at (3,0) {$\Omega K_{t,\mathbb{C}}$};
\node (Ru) at (3,2) {$\Omega K_{t}$};

\draw [->](Ll) to node [above]{\scriptsize $t_{h}^{\natural}$}(Rl); 
\draw [->](Ru) to (Ll);
\draw [->](Ru) to (Rl);
\end{tikzpicture}
\end{center}
gives an infinite loop map 
\[e_{h}^{\natural}:K_{a}\mathbb{C}\rightarrow F_{t,\mathbb{C/Z}}\] which makes the diagram \eqref{tecommutativediag} commute. Combining the diagram \eqref{tecommutativediag} with the fact that $\pi_{\ast}(K_{t})=0$, for $\ast$ is odd---hence $\pi_{\ast}(K^{rel}\mathbb{C})\rightarrow \pi_{\ast}(K_{a}\mathbb{C})$ is onto, we obtain:  
\[e_{h,\ast}^{\natural}=e_{\ast}:\pi_{\ast}(K_{a}\mathbb{C})\rightarrow \pi_{\ast}(F_{t,\mathbb{C/Z}}).\] Notice, for $\ast$ is even, $\pi_{\ast}(F_{t,\mathbb{C/Z}})=0$.

As for the uniqueness, we assume there is another lifting $e_{h}^{\flat}$ such that 
\[e_{h,\ast}^{\flat}=e_{\ast}=e_{h,\ast}^{\natural}.\]
In view of Lemma \ref{edeterminest}, we may assume $t_{h}^{\flat}$ is the induced infinite loop map from $K^{rel}\mathbb{C}$ to $\Omega K_{t,\mathbb{C}}$. Let $p_{\ast}$ be the homomorphism 
\[\pi_{\ast}(\Omega K_{t,\mathbb{C}})\rightarrow \pi_{\ast}(F_{t,\mathbb{C/Z}}).\]
Then the diagram \eqref{tecommutativediag} implies  
\begin{equation}\label{EqforEflatandEnatural}
p_{\ast}\circ t^{\flat}_{h,\ast}=e^{\flat}_{\ast}\circ \pi_{\ast}=e^{\natural}_{\ast}\circ \pi_{\ast}=p_{\ast}\circ t^{\natural}_{h,\ast}.
\end{equation}   
Now observe there is an exact sequence 
\begin{multline*}
0\rightarrow \operatorname{Hom}(\pi_{\ast}(K^{rel}\mathbb{C}),\pi_{\ast}(\Omega K_{t}))\rightarrow \operatorname{Hom}(\pi_{\ast}(K^{rel}\mathbb{C}),\pi_{\ast}(\Omega K_{t,\mathbb{C}}))\\
\rightarrow \operatorname{Hom}(\pi_{\ast}(K^{rel}\mathbb{C}),\pi_{\ast}(F_{t,\mathbb{C/Z}}))
\end{multline*}
given by the short exact sequence
\[0\rightarrow \pi_{\ast}(\Omega K_{t})\rightarrow \pi_{\ast}(\Omega K_{t,\mathbb{C}})\rightarrow \pi_{\ast}(F_{t,\mathbb{C/Z}})\rightarrow 0.\]
Since $\operatorname{Hom}(\pi_{\ast}(K^{rel}\mathbb{C}),\pi_{\ast}(\Omega K_{t}))=0$, the homomorphism 
\[\operatorname{Hom}(\pi_{\ast}(K^{rel}\mathbb{C}),\pi_{\ast}(\Omega K_{t,\mathbb{C}}))\rightarrow \operatorname{Hom}(\pi_{\ast}(K^{rel}\mathbb{C}),\pi_{\ast}(F_{t,\mathbb{C/Z}}))\]
is actually injective. Hence, the maps $t^{\flat}_{h,\ast}=t^{\natural}_{h,\ast}$, in view of the equality \eqref{EqforEflatandEnatural}, and $t_{h}^{\flat}$ and $t_{h}^{\natural}$ are homotopic in $\operatorname{Ho}(\mathcal{P})$ (see \cite[Lemma $2.2.7$]{Wang1}). Applying Lemma \ref{thdetermineseh}, we obtain the maps $e^{\flat}_{h}$ and $e^{\natural}_{h}$ differ only by a phantom map.

The assertion \eqref{Eq:Torenande} has been shown in Corollary\ref{Torsionpartarethesame}, whereas the equality \eqref{Eq:Boanden} follows from the definition of $t^{\natural}_{h}$ and the fact that $\operatorname{Bo}$ is also an infinite loop map as we have 
\[\operatorname{ch}_{\otimes\mathbb{R},\ast}\circ \operatorname{Im}_{\ast}\circ t^{\natural}_{h,\ast}=\operatorname{ch}_{\otimes\mathbb{R},\ast}\circ \operatorname{Im}_{\ast} \circ \operatorname{ch}^{rel}_{\ast}=\operatorname{Bo}_{\ast}:\pi_{\ast}(K^{rel}\mathbb{C})\rightarrow \pi_{\ast}(H^{odd}\mathbb{R}).\]
\end{proof}
The theorem above gives strong evidence for the conjecture that $e$ can be lifted to a map in the stable homotopy category. In fact, if one can show $\operatorname{ch}^{rel}$ is an infinite loop map and hence $\operatorname{ch}^{rel}= t_{h}^{\natural}\in \operatorname{Ho}(\mathcal{P})$, then $e_{h}^{\natural}$ and $e$ are weakly homotopic as all free elements in $[L,K_{a}\mathbb{C}]$ come from $[L,K^{rel}\mathbb{C}]$ when $L$ is a finite $CW$-complex.

%%%%%%%%%%%%%%%%%%%%%%%%

\section{Index theorems for flat vector bundles}
Combining Theorem \ref{enatural} with the commutative diagram \eqref{Intro:DWWindex}, which is given by the $\operatorname{DWW}$ index theorem, we have the following refined 
$\operatorname{BL}$ index theorem:

\begin{theorem}\label{IndexThm1}
Given a compact smooth fiber bundle $E\rightarrow B$, then the diagram  
\begin{center}
\begin{equation} 
\begin{tikzpicture}[baseline=(current  bounding  box.center)]
\node(Lu) at (0,2) {$\tilde{K}(E,\mathbb{C})$};
\node(Ll) at (0,0) {$\tilde{K}(B,\mathbb{C})$}; 
\node(Ru) at (4,2) {$[E,F_{t,\mathbb{C/Z}}]$};
\node(Rl) at (4,0) {$[B,F_{t,\mathbb{C/Z}}]$};

\draw[->] (Lu) to node [right]{\scriptsize $\pi^{!}$}(Ll);
\draw[->] (Ru) to node [right]{\scriptsize $\operatorname{tr}_{\operatorname{BG}}^{\ast}$}(Rl);
\draw[->] (Lu) to node [above]{\scriptsize $\bar{e}^{\natural}$}(Ru);
\draw[->] (Ll) to node [above]{\scriptsize $\bar{e}^{\natural}$}(Rl);  
\end{tikzpicture}
\end{equation}
\end{center}
commutes, where $\bar{e}^{\natural}$ is the composition
\[\tilde{K}(-,\mathbb{C})\rightarrow [-,K_{a}\mathbb{C}]\xrightarrow{e^{\natural}_{\ast}} [-,F_{t,\mathbb{C/Z}}].\]
\end{theorem}

Theorem \ref{IndexThm1} implies the following index theorem in terms of $\bar{e}_{\operatorname{APS}}$:
\begin{theorem}
Let $E\rightarrow B$ be a smooth compact fiber bundle with the fundamental group $\pi_{1}(E,\ast)$ finite, for every point $\ast\in E$, then the diagram   
\begin{center}
\begin{equation} 
\begin{tikzpicture}[baseline=(current  bounding  box.center)]
\node(Lu) at (0,2) {$\tilde{K}(E,\mathbb{C})$};
\node(Ll) at (0,0) {$\tilde{K}(B,\mathbb{C})$}; 
\node(Ru) at (4,2) {$[E,F_{t,\mathbb{C/Z}}]$};
\node(Rl) at (4,0) {$[B,F_{t,\mathbb{C/Z}}]$};

\draw[->] (Lu) to node [right]{\scriptsize $\pi^{!}$}(Ll);
\draw[->] (Ru) to node [right]{\scriptsize $\operatorname{tr}_{\operatorname{BG}}^{\ast}$}(Rl);
\draw[->] (Lu) to node [above]{\scriptsize $\bar{e}_{\operatorname{APS}}$}(Ru);
\draw[->] (Ll) to node [above]{\scriptsize $\bar{e}_{\operatorname{APS}}$}(Rl);  
\end{tikzpicture}
\end{equation}
\end{center} 
commutes. 
\end{theorem}
\begin{proof}
The assumption that the fundamental group $\pi_{1}(E,\ast)$ is finite, for every point $\ast\in E$, implies the fundamental group $\pi_{1}(B,\ast)$ is also finite, for every point $\ast\in B$, and it is also known that $[BG,K(\mathbb{Q},i)]$ is trivial, for any finite group $G$ and every $i\in\mathbb{N}$, where $BG$ is the classifying space of $G$ (see \cite[Corollary $4.3$]{Web}). Thus, given any representation $G\rightarrow \operatorname{GL}(\mathbb{C})$, the induced homomorphism $[X,BG]\rightarrow [X,K_{a}\mathbb{C}]$ always factors through $\operatorname{Tor}[X,K_{a}\mathbb{C}]$, where $X$ is a topological space, as there are isomorphisms: 
\[[BG,K_{a}\mathbb{C}]\otimes\mathbb{Q}\simeq [BG,K_{a}\mathbb{C}_{\mathbb{Q}}]\simeq [BG,\prod_{i\in\mathbb{N}}K(\pi_{i}(K_{a}\mathbb{C})\otimes\mathbb{Q},i)]=0.\]

The assertion then follows from the commutative diagram below:
\begin{center}
\begin{equation} 
\begin{tikzpicture}[baseline=(current  bounding  box.center)]
\node(Lu) at (0,2) {$\tilde{K}(E,\mathbb{C})$};
\node(Ll) at (0,0) {$\tilde{K}(B,\mathbb{C})$}; 
\node(Mu) at (3,2) {$\operatorname{Tor}[E,K_{a}\mathbb{C}]$}; 
\node(Ml) at (3,0) {$\operatorname{Tor}[B,K_{a}\mathbb{C}]$};
\node(Ru) at (7,2) {$[E,F_{t,\mathbb{C/Z}}]$};
\node(Rl) at (7,0) {$[B,F_{t,\mathbb{C/Z}}]$};

\draw[->] (Lu) to node [right]{\scriptsize $\pi^{!}$}(Ll);
\draw[->] (Lu) to   (Mu);
\draw[->] (Ll) to   (Ml);
\draw[->] (Lu) to [out=30,in=150]node [above]{\scriptsize $\bar{e}_{\operatorname{APS}}$}  (Ru);
\draw[->] (Ll) to [out=-30,in=-150]node [above]{\scriptsize $\bar{e}_{\operatorname{APS}}$}  (Rl);
\draw[->] (Mu) to node [right]{\scriptsize $\operatorname{tr}_{\operatorname{BG}}^{\ast}$} (Ml); 
\draw[->] (Ru) to node [right]{\scriptsize $\operatorname{tr}_{\operatorname{BG}}^{\ast}$}(Rl);
\draw[->] (Mu) to node [above]{\scriptsize $e_{\ast}= e^{\natural}_{h.\ast}$}(Ru);
\draw[->] (Ml) to node [above]{\scriptsize $e_{\ast}= e^{\natural}_{h,\ast}$}(Rl);  
\end{tikzpicture}
\end{equation}
\end{center} 

\end{proof}

\section{The Adams $e$-invariant}
%%%%%%%%The method must change! 
In this section, we explain how the $e$-invariant \cite[Section $4.1$]{Wang1} generalizes the Adams $e$-invariant. Recall the Adams $e$-invariant can be obtained as the lifting of the following diagram---there exists only one lifting as we have $[B\Sigma_{\infty},K(\mathbb{Q},i)]=0$, for every $i\in\mathbb{N}$ \cite[Corollary $4.3$]{Web}: 
\begin{center}
\begin{tikzpicture}
\node(Ll) at (0,0) {$B\Sigma_{\infty}$};
\node(Ml) at (2,0) {$K_{t}$};
\node(Mu) at (2,2) {$F_{t,\mathbb{Q/Z}}$};
\node(Ru) at (4,2) {$F_{t,\mathbb{C/Z}}$};

\path[->, font=\scriptsize,>=angle 90]

(Ll) edge (Ml)
(Ll) edge node [xshift=-3ex,yshift=1ex]{$\bar{e}_{\operatorname{Adams}}$}(Mu)
(Mu) edge (Ml)
(Ru) edge (Ml)
(Mu) edge (Ru);

\end{tikzpicture}
\end{center}
Applying the universal property of the plus construction, one gets a map
\[e_{\operatorname{Adams}}:B\Sigma_{\infty}^{+}\rightarrow F_{t,\mathbb{Q/Z}}\]
whose induced homomorphism
\[e_{\operatorname{Adams},\ast}:\pi_{\ast}(B\Sigma_{\infty}^{+})\rightarrow \pi_{\ast}(F_{t,\mathbb{Q/Z}})\]
gives the Adams $e$-invariant up to sign \cite{Qu2}.

Because $[B\Sigma_{\infty},\Omega K_{t,\mathbb{C}}]$ is trivial, the following diagram must commute:
\begin{center} 
\begin{tikzpicture}
\node(Ll) at (0,0){$B\Sigma_{\infty}^{+}$}; 
\node(Ml) at (2.5,0){$K_{a}\mathbb{C}$};
\node(Mu) at (2.5,2){$F_{t,\mathbb{Q/Z}}$};
\node(Rl) at (5,0){$K_{t}$};
\node(Ru) at (5,2){$F_{t,\mathbb{C/Z}}$};
 
\path[->, font=\scriptsize,>=angle 90] 

(Ll) edge node [xshift=-2.2ex,yshift=1ex] {$e_{\operatorname{Adams}}$} (Mu)
(Mu) edge (Ru)
(Ll) edge node [above]{$\iota$}(Ml)
(Ml) edge (Rl)
(Ml) edge node [xshift=-2ex,yshift=.5ex]{$e_{h}^{\mathbb{C}}$}(Ru)
(Ru) edge (Rl);

\end{tikzpicture}
\end{center} 
where $e_{h}^{\mathbb{C}}$ is any lifting of the comparison map $K_{a}\mathbb{C}\rightarrow K_{t}$. In particular, this shows the two homomorphisms  
\begin{align*}
[L,B\Sigma^{+}_{\infty}]&\xrightarrow{\iota_{\ast}} [L,K_{a}\mathbb{C}]\xrightarrow{e^{\mathbb{C}}_{h,\ast}} [L,F_{t,\mathbb{C/Z}}];\\ 
[L,B\Sigma^{+}_{\infty}]&\xrightarrow{e_{\operatorname{Adams},\ast}} [L,F_{t,\mathbb{Q/Z}}]\rightarrow [L,F_{t,\mathbb{C/Z}}],\\ 
\end{align*}
are the same, for any finite $CW$-complex $L$. Furthermore, by Corollary \ref{sameisoLiftings}, we know $e^{\mathbb{C}}_{h,\ast}=e_{\ast}$ in this case. Thus, we have proved the following theorem:

\begin{theorem}
\[(e\circ\iota)_{\ast}=e_{\operatorname{Adams},\ast}:\pi_{\ast}(B\Sigma_{\infty}^{+})\rightarrow \pi_{\ast}(F_{t,\mathbb{C/Z}}).\] 
\end{theorem} 
  
%%%Try to avoid CW-spectra
%%%Using homotopy theoretic way to prove it generalizes teh Adams e-invairant

%%%State some possible ways to further investigate the problem
%--The aims is to show tch is an infinite loop space
%Consider the summand
%Probably its real part contains no information about torsion-free 
%Or maybe when its holonomy group is trivial, the real part contains no information about torsion free part.

%%%%%%%%%%%%%%%%%%%%%%%%%%%%%%%%%%%%%%%%%%%%%%%%%%%%%%%%%%%%%%%%%%%%%%%%%%%%%%%%%%%
%%%%%%%%%%%%%%%%%%%%%%%%%%%%%%%%%%%%%%%%%%%%%%%%%%%%%%%%%%%%%%%%%%%%%%%%%%%%%%%%%%%%%%%
%%%%%%%%%%%%%%%%%%%%%%%%%%%%%%%%%%%%%%%%%%%%%%%%%%%%%%%%%%%%%%%%%%%%%%%%%%%%%%%%%%%%%%%%%
%%%%%%%%%%%%%%%%%%%%%%%%%%%%%%%%%%%%%%%%%%%%%%%%%%%%%%%%%%%%%%%%%%%%%%%%%%%%%%%%%%%%%%%%%
%%%%%%%%%%%%%%%%%%%%%%%%%%%%%%%%%%%%%%%%%%%%%%%%%%%%%%%%%%%%%%%%%%%%%%%%%%%%%%%%%%%%%%%%%%
%%%%%%%%%%%%%%%%%%%%%%%%%%%%%%%%Main Body begins$$$$$$$$$$$$$$$$$$$$$$$$$$$$$$$$$$$$$$$

\newpage
{\center \textbf{Acknowledgment}

The results of the present paper improve some theorems of \cite{Wang3}. The author would like to thank his supervisor Sebastian Goette for bring to his attention this interesting research project. The author is also grateful to Ulrich Bunke and Wolfgang Steimle for their expert advice and insightful comments on the author's work. The author would also like to express his sincere thanks to Jørgen Olsen Lye for reading part of the early draft carefully. The project is financially supported by the DFG Graduiertenkolleg 1821 "Cohomological Methods in Geometry".

}
%%%%%%%%%%%%%%%%%%%%%%%%%%%%%%%%%%%%%%%%%%%%%%%%%%%%%%%%%%%%%%%%%%%%%%%%%%%%%%%%%%%%%%%%%%%%%%%%%%%%%%%%%%%%%%%%%%
%%%%%%%%%%%%%%%%%%%%%%%%%%%%%%%%%%%%%%%%%%%%%%%%%%%%%%%%%%%%%%%%%%%%%%%%%%%%%%%%%%%%%%%%%%%%%%%
%%%%%%%%%%%%%%%%%%%%%%%%%%%%%%%%%%%%%%%%%%%%%%%%%%%%%%%%%%%%%%%%%%%%%%%%%%%%%%%%%%%%%%%%%%%%%%%
  
\phantomsection 
\addcontentsline{toc}{section}{Bibliography}

\bibliographystyle{alpha}
 
\bibliography{Reference}

\end{document}